\documentclass[12pt]{amsart}
\usepackage{amsmath}
\usepackage{amsfonts}
\usepackage{amssymb}
\usepackage{amscd}
\catcode `\@=11
\def\numberbysection{\@addtoreset{equation}{section}
         \renewcommand{\theequation}{\thesection.\arabic{equation}}}
\numberbysection
\def\subsubsection{\@startsection{subsubsection}{3}%
  \normalparindent{.5\linespacing\@plus.7\linespacing}{-.5em}%
  {\normalfont\bfseries}}
\def\a{\alpha}
\def\b{\beta}
\def\g{\gamma}

\def\eps{\epsilon}
\def\l{\lambda}
\def\s{\sigma}
\def\t{\tau}

\def\nn{\nonumber}
\def\la{\langle}
\def\ra{\rangle}
\def\Z2{\mathbb {Z}/2\mathbb{Z}}

\def\Str{\mathrm{Str}}
\def\Hom{\mathrm{Hom}}

\begin{document}

\title{The algebra of discrete torsion}

\author[Ralph M. Kaufmann]{Ralph M. Kaufmann$^*$\\
University of Southern California, Los Angeles, USA\\
and Max--Planck Institut f\"ur Mathematik, Bonn, Germany}

\thanks{${}^*$ Partially supported by NSF grant \#0070681}
\email{kaufmann@math.usc.edu}

\address{University of Southern California, Los Angeles, USA and
 Max--Planck Institut f\"ur Mathematik, Bonn, Germany}

\date{\today}
\begin{abstract}
We analyze the algebraic structures of $G$--Frobenius algebras
which are the algebras associated to global group quotient
objects. Here $G$ is any finite group. These algebras turn out to
be modules over the Drinfeld double of the group ring $k[G]$. We
furthermore prove that discrete torsion is a universal group
action of $H^2(G,k^*)$ on $G$--Frobenius algebras by isomorphisms
of the underlying linear structure. These morphisms are realized
explicitly by taking the tensor product with twisted group rings.
This gives an algebraic realization of discrete torsion and allows
for a treatment analogous to the theory of projective
representations of groups, group extensions and twisted group ring
modules. Lastly, we identify another set of discrete universal
transformations among $G$--Frobenius algebras pertaining to their
super--structure and classified by $\Hom(G,\Z2)$.
\end{abstract}
\maketitle

\section*{Introduction}

$G$--Frobenius algebras were introduced by the author in [K1] to
explain the algebraic structure replacing Frobenius algebras when one
is dealing with global group quotients in theories such as Cohomology,
K--theory, local rings of singularities etc.. The main feature
is that they are group graded non--commutative algebras with a
group action and a controlled non--commutativity.
Previously, one had only regarded commutative algebra structures
on the invariants under the group action of the $G$--module,
which leads to a direct sum decomposition indexed
by conjugacy classes. The larger non--commutative
 algebras however reflect the geometry and general properties much better.
For instance in the present treatise it is the full $G$--Frobenius
structure which allows us to decode discrete torsion in all its
aspects. Our $G$--Frobenius algebras have been proven to be
exactly the right structure to describe the cohomology of global
quotients [FG] and the geometry of symmetric products [K2,LS]
whose structure is closely related to Hilbert schemes [LS].
Furthermore they are suited to describe mirror phenomena for
Frobenius structures having their origin in singularity theory
[K1,K3].  One fundamentally new feature of these algebras is that
given any one of them one expects the existence of several cousins
of it, which are all related to each other by discrete torsion.
There are several definitions of discrete torsion
[V,VW,D,A,R,S,K1,K2]. The purpose of this note is to show that all
definitions can be understood in terms of the fact that there is a
universal group action of $H^2(G,k^*)$ on $G$--Frobenius algebras,
by isomorphisms of the underlying linear structure, by taking
tensor products with twisted group rings.

The algebraic structures of a $G$--Frobenius algebra are as follows:
it is naturally a left $k[G]$--module
algebra as well as a right $k[G]$--co--module algebra. Moreover it
satisfies the Yetter--Drinfeld (YD) condition for bi--modules and
is thus a module over $D(k[G])$, the Drinfeld double of $k[G]$.

We show that given a $G$--Frobenius algebra there is an
action of $H^2(G,k^*)$ on it by tensoring with the respective twisted
group ring. This tensor product exhibits all desiderata of discrete
torsion. In fact this action is an action of
universal projective re-scalings.
Moreover the action transforms the partition sums in the right way and
defines a bi--character (usually denoted by $\eps(g,h)$) on commuting
elements with the expected properties. It not only satisfies the
necessary algebraic relations, but it also appears
naturally as  factors in front of the summands of the partition function.

The reason for the beautiful picture is the inclusion of twisted sectors
for all group elements and not just for conjugacy classes, which
has been discussed in detail in [K1,K2]. This also allows to understand
the bi--character in terms of a 2 cocycle $\a \in H^2(G,k^*)$.
The bi--character
$\eps$ is then derived from $\a$ and is defined for all elements.
It can actually be shown to be a one cocycle in $\eps \in H^1(G,k^*[G])$
where $k^*[G]\subset k[G]$ are sums with invertible coefficients
and the $G$--module structure is
given by the adjoint action.

We furthermore comment on the generic super--structures one can
impose on a given $G$--Frobenius algebra and realize
that these are also given by tensor product with
superized versions of $k[G]$. These are a second type of
discrete deformation, which is actually different from the one
of discrete torsion.

This action of discrete torsion using twisted group rings
should be viewed as a statement which is analogous to
the relations between projective representations of a group,
modules over the twisted group algebra and extensions of the group.

In fact we are able to establish the counterpart
to the classical picture in the world of
$G$--Frobenius algebras.

The paper is organized as follows:

In the first paragraph, we recall the notion of $G$--Frobenius algebras,
and show that group rings and twisted group rings are
$G$--Frobenius algebras.

In the second paragraph, we give the algebraic properties of $G$--Frobenius
algebras. The main Theorem is that a $G$--Frobenius algebra
has the natural algebraic structures stated above.
We also characterize $G$--Frobenius algebras which are
Galois over their identity sector as $k[G]$--comodule algebras.

The third section contains the identification of discrete
torsion as an action of $H^{2}(G,k^*)$ on $G$--Frobenius algebras.
This is done by analyzing universal twists of $G$--Frobenius algebras and
showing that the universal twists are in 1--1 correspondence with $H^2(G,k^*)$.
Furthermore we show that these twists can be realized by tensoring
with the respective twisted group ring.

We furthermore also study the generic super--structures ($\Z2$) which
one can impose on a given $G$--Frobenius algebra.

In the fourth paragraph we develop our understanding
of these facts by introducing a theory in analogy with
 projective representations of a group and their relation to
modules over the twisted group algebra and extensions of the group.
Here the final result is that given any Abelian group $H$ and a
cocycle $[\a'] \in H^2(G,H)$ then for any central extension $G^{\a}$ of
$G$ by $H$ with class $[\a']$ and any $G$--Frobenius algebra $A$
there is a natural $G^{\a'}$--Frobenius algebra $A^{\a'}$ to which
the Frobenius algebra $A_{\a}$ ($A$ twisted by $\a$) can be lifted.
Here $[\a]\in H^2(G,k^*)$ is the image under the transgression map
associated to $[\a']$ of a $\chi \in \Hom(H,k^*)$.
Vice versa, the above $\chi$ gives a push down map, which maps
$A^{\a}$ onto $A_{\a}$.
Furthermore there is a universal setup of this kind if there is a
representation group for $G$.

\section*{Acknowledgements}
I would like to thank the Max--Planck--Institut for its kind hospitality
and also gratefully  acknowledge the support from the NSF.
Special thanks goes to Takashi Kimura, since it was a discussion with
him that was the initial spark for this paper.
I would also like to thank Susan Montgomery whose questions about
the algebraic setting of orbifolding are the reason for the
existence of the second section, which also owes a great deal to her
book on ``Hopf algebras and their actions on rings''.

\section{$G$--Frobenius algebras}
\label{orb}

We fix a finite group $G$ and denote its unit element by $e$. We
furthermore fix a ground field $k$ of characteristic zero for
simplicity. With the usual precautions the characteristic
of the field does not play an important role and furthermore
the group really only needs to be completely disconnected.

\subsection{Definition}
 A {\em G--twisted Frobenius algebra} ---or $G$--Frobenius algebra for short---
over
a field  $k$ of characteristic 0 is
$<G,A,\circ,1,\eta,\varphi,\chi>$, where

\begin{tabular}{ll}
$G$&finite group\\
$A$&finite dim $G$-graded $k$--vector space \\
&$A=\oplus_{g \in G}A_{g}$\\
&$A_{e}$ is called the untwisted sector and \\
&the $A_{g}$ for $g \neq
e$ are called the twisted sectors.\\
$\circ$&a multiplication on $A$ which respects the grading:\\
&$\circ:A_g \otimes A_h \rightarrow A_{gh}$\\
$1$&a fixed element in $A_{e}$--the unit\\
$\eta$&non-degenerate bilinear form\\
&which respects grading i.e. $g|_{A_{g}\otimes A_{h}}=0$ unless
$gh=e$.\\
\end{tabular}

\begin{tabular}{ll}
$\varphi$&an action  of $G$ on $A$
(which will be  by algebra automorphisms), \\
&$\varphi\in \mathrm{Hom}(G,\mathrm{Aut}(A))$, s.t.\
$\varphi_{g}(A_{h})\subset A_{ghg^{-1}}$\\
$\chi$&a character $\chi \in \mathrm {Hom}(G,k^{*})$ \\

\end{tabular}

\vskip 0.3cm

\noindent Satisfying the following axioms:

\noindent{\sc Notation:} We use a subscript on an element of $A$ to signify that it has homogeneous group
degree  --e.g.\ $a_g$ means $a_g \in A_g$-- and we write $\varphi_{g}:= \varphi(g)$ and $\chi_{g}:= \chi(g)$.

\begin{itemize}

\item[a)] {\em Associativity}

$(a_{g}\circ a_{h}) \circ a_{k} =a_{g}\circ (a_{h} \circ a_{k})$
\item[b)] {\em Twisted commutativity}

$a_{g}\circ a_{h} = \varphi_{g}(a_{h})\circ a_{g}$
\item[c)]
{\em $G$ Invariant Unit}:

$1 \circ a_{g} = a_{g}\circ 1 = a_g$

and

$\varphi_g(1)=1$
\item[d)]
{\em Invariance of the metric}:

$\eta(a_{g},a_{h}\circ a_{k}) = \eta(a_{g}\circ a_{h},a_{k})$

\item[i)]
{\em Projective self--invariance of the twisted sectors}

$\varphi_{g}|A_{g}=\chi_{g}^{-1}id$

\item[ii)]
{\em $G$--Invariance of the multiplication}

$\varphi_{k}(a_{g}\circ a_{h}) = \varphi_{k}(a_{g})\circ  \varphi_{k}(a_{h})$

\item[iii)]

{\em Projective $G$--invariance of the metric}

$\varphi_{g}^{*}(\eta) = \chi_{g}^{-2}\eta$

\item[iv)]
{\em Projective trace axiom}

$\forall c \in A_{[g,h]}$ and $l_c$ left multiplication by $c$:

$\chi_{h}\mathrm {Tr} (l_c  \varphi_{h}|_{A_{g}})=
\chi_{g^{-1}}\mathrm  {Tr}(  \varphi_{g^{-1}} l_c|_{A_{h}})$
\end{itemize}

\subsubsection{Special $G$--Frobenius algebras}
We briefly review special $G$--Frobenius algebras. For details
see [K1,K2].

We call a $G$-Frobenius algebra special if all $A_g$ are cyclic
$A_e$ modules via the multiplication $A_e \otimes A_g \rightarrow A_g$.
Fixing a cyclic generator $1_g \in A_g$ the
algebra is completely characterized by two compatible
cocycles, namely $\g \in \bar Z^2(G,A_e)$ and $\varphi\in Z^1(G,k^*[G])$
 where $\bar Z$ are graded cocyles (see [K1]) and
$k^*[G]$ is the group ring restricted to invertible coefficients
with $G$--module structure induced by the adjoint action:
$$
\phi(g)\cdot(\sum \mu_h h)= \sum_h \mu_h ghg^{-1}
$$
We set $\varphi(g) = \sum_h \varphi_{g,h} ghg^{-1}$ and
$\g_{g,h} = \g(g,h)$.

The multiplication and $G$--action are determined by
$$
1_g 1_h = \g_{g,h} 1_{gh} \quad \varphi_{g}(1_{h}) = \varphi_{g,h}
1_{ghg^{-1}}
$$
There are two compatibility equations:
\begin{equation}
\label{grpcompat}
\varphi_{g,h}\g_{ghg^{-1},g} = \g_{g,h}
\end{equation}
and
\begin{equation}
\label{algaut}
\varphi_{k,g} \varphi_{k,h} \g_{kgk^{-1},khk^{-1}}
= \varphi_{k} (\g_{g,h}) \varphi_{k,gh}
\end{equation}

Notice that if $\g_{g,h}$ is non--zero i.e. $A_gA_h \neq 0$ then
(\ref{grpcompat})
determines $\varphi_{g,h}$. We also would like to remark that
(\ref{algaut}) is automatically satisfied if $A_gA_hA_k \neq 0$
(cf.\ [K1]).

\subsection{The group ring $k[G]$}
Let $k[G]$ denote the group ring of $G$.
\subsubsection{The Hopf structure of $k[G]$} Recall that $k[G]$ is a Hopf algebra with
the natural multiplication, the co--multiplication induced by
$\Delta(g) = g \otimes g$, co--unit $\eps(g)=1$ and antipode $S(g)
= g^{-1}$.

\subsubsection{The $G$--Frobenius structure of $k[G]$}
When considering $k[G]$ as a $G$--Frobenius algebra we will
consider $k[G]$ as a left $k[G]$--module with respect to
conjugation, i.e. the map $k[G] \otimes k[G] \rightarrow k[G]$
given by

$$
\sum_g \nu_g g \otimes \sum_h \mu_h h \mapsto \sum_{g,h}\nu_h
\mu_g ghg^{-1}
$$

The other structures are the naturally $G$ graded natural
multiplication on $k[G]$ with the unit $e$, the metric
$\eta(g,h)=\delta_{gh,e}$ and $\chi_g \equiv 1$. It is trivial to
check all axioms.

If we were to choose a grading $\tilde{}\in \Hom(G,\Z2)$ then
$\chi_g= (-1)^{\tilde g}$ and $\varphi_{g}(h) =(-1)^{\tilde g
\tilde h}$.

\subsection{The twisted group ring $k^{\a}[G]$}

Recall that given an element $\a \in Z^2(G,k^*)$
one defines the twisted group ring
$k^{\a}[G]$ to be given by the same linear structure with multiplication
given by the linear extension of

\begin{equation}
g\otimes h \mapsto \a(g,h) gh
\end{equation}
with $1$ remaining the unit element.
To avoid confusion we will denote elements of $k^{\a}[G]$ by
$\hat g$ and the multiplication with $\cdot$
Thus
$$
\hat g \cdot \hat h = \a(g,h) \widehat{gh}
$$

For $\a$ the following equations hold:
\begin{equation}
\a (g,e) = \a(e,g), \qquad
\a(g,g^{-1})=\a(g^{-1},g)
\end{equation}
Furthermore
$$
\hat{g}^{-1}= \frac{1}{\a(g,g^{-1})}\hat{g^{-1}}
$$

\subsubsection{Remark}
Given a two co--cycle $\a$ and
possibly extending the field by square roots
we can find a co--cycle $\tilde \a$  in the same
cohomology class which also satisfies
\begin{equation}
\tilde \a(g,g^{-1})=1
\end{equation}
If one wishes to consider $\mathbb{C}$ as a ground field, one can
work with such cocycles.
\subsubsection{Lemma}
\label{adjointaction} Set $\eps(g,h)=\frac{\a(g,h)}{\a(ghg^{-1},g)}$,
then the left adjoint action of $k^{\a}[G]$ on $k^{\a}[G]$ is
given by
\begin{equation}
g\otimes h \stackrel{ad}{\mapsto} \eps(g,h) \widehat{ghg^{-1}}
\end{equation}

{\bf Proof.} By the definition of multiplication in $k^{\a}[G]$
$$
\hat g\cdot \hat h\cdot\hat {g}^{-1} =
\frac{\a(g,h)\a(gh,g^{-1})}{\a({g,g^{-1})}}\widehat{ghg^{-1}}
$$
Now
 by associativity
$$
\a(gh,g^{-1})\a(ghg^{-1},g)=\a(gh,e)\a(g^{-1},g)=\a(g,g^{-1}).
$$
so
$$
\frac{\a(g,h)\a(gh,g^{-1})}{\a(g,g^{-1})}=\frac{\a(g,h)}{\a(ghg^{-1},g)}
$$

\subsubsection{The $G$--Frobenius Algebra structure of $k^{\a}[G]$}

Recall from [K1,K2] the following structures which turn
$k^{\a}[G]$ into a special $G$--Frobenius algebra:

\begin{eqnarray}
 \g_{g,h}&=&\a(g,h)\nn\\
\varphi_{g,h}&=& (-1)^{\tilde g\tilde h}
\frac{\a(g,h)}{\a(ghg^{-1},g)}=:\eps(g,h)\nn\\
\eta(\hat g,\hat{g^{-1}}) &=&\a(g,g^{-1})\nn\\
\chi_g&=& (-1)^{\tilde g}
\end{eqnarray}

Here the second line induces the third via
$$
\hat g\cdot \hat h = \a(g,h)\hat{gh}, \quad  \widehat{ghg^{-1}} \cdot
\hat g = \a(ghg^{-1},g)\hat{gh}
$$

We recall that if $k^*$ is two divisible we could scale s.t.\
$\eta(g,g^{-1})=1$ and $\eps$ would indeed yield the adjoint action.
The last equation follows from the special case of the trace axiom
since the dimension of all sectors is one.

It is an exercise to check all axioms. All compatibility equations follow
automatically, since $\a(g,h) \neq 0$.
The only axiom which is not straightforward is the trace axiom, but see
[K2] for a proof.

\subsubsection{Remark} By the general theory (see above),
$\eps \in H^1(G,k^*[G])$ where
$k^*[G]$ is the group ring restricted to invertible coefficients
with $G$--module structure induced by the adjoint action:
$$
\phi(g)\cdot(\sum \mu_h h)= \sum_h \mu_h ghg^{-1}
$$

\subsubsection{Relations}
The $\eps(g,h)$ satisfy the equations:
\begin{eqnarray}
\eps(g,e)&=&\eps(g,g)=1\\
\label{group} \eps(g_1g_2,h)&=&
\eps(g_1,g_2hg_2^{-1})\eps(g_2,h)\\
\label{invert} \eps(g,h)^{-1}&=&\eps(g^{-1},ghg^{-1})\\
\label{compat} \eps(k,gh) &=&
\eps(k,g)\eps(k,h)\frac{\a(kgk^{-1},khk^{-1})}{\a(g,h)}
\end{eqnarray}

where (\ref{group}) is the statement that $\varphi \in \Hom(G,
Aut(A))$, (\ref{invert}) is a consequence of (\ref{group}) and
(\ref{compat}) is the compatibility equation which also ensures
the invariance of the metric.

Furthermore the trace axiom holds [K2] which is equivalent to the
equation

\begin{equation}
\a([g,h],hgh^{-1})\eps(h,g)= \a([g,h],h)\eps(g^{-1},ghg^{-1})
\end{equation}
or
\begin{equation}
\eps(h,g)=\eps(g^{-1},ghg^{-1})
\frac{\a([g,h],h)}{\a([g,h],hgh^{-1})}
\end{equation}

 In the case that the group elements in the equations commute we
obtain the famous conditions of discrete torsion which make $\eps$
into a bi--character on commuting elements.

For {\em commuting elements}:

\begin{eqnarray}
\eps(g,e)&=&\eps(g,g)=1\nn\\
\eps(g_1g_2,h)&=& \eps(g_1,h)\eps(g_2,h)\nn\\
\eps(g,h)^{-1}&=&\eps(g^{-1},h)\nn\\
 \eps(g,h) &=&\eps(h^{-1},g)=\eps(h,g)^{-1}\nn\\
  \eps(h,g_1g_2) &=& \eps(h,g_1) \eps(h,g_2)
\end{eqnarray}
where the last equation is now a consequence of the second and the
fourth and the third equation follows from the second.

\subsubsection{Fact} One can show [K2] that the twisted group
algebras $k^{\a}[G]$ are the only $G$--Frobenius algebras with
the property that all $A_g$ are one--dimensional.
To be completely precise there is an additional
freedom of choosing a super (i.e.\ $\Z2$) structure
determined by a homomorphism $\s \in \Hom(G,\Z2)$
(see [K1] and  \ref{super} below).

\subsubsection{Geometry of $k^{\a}[G]$} From the point of view
of Jacobian Frobenius algebras [K1] it is natural to
say that $k[G]$ is the Frobenius algebra naturally associated
to $point/G$. The existence of the twisted algebras suggests
that there are several equivalent ways of taking the group quotient.
This is made precise by Theorem \ref{dt} below.

\section{Algebraic structures of a $G$--Frobenius algebra}

We fix a $G$--Frobenius algebra $\la G,A,\circ,1,\eta,\varphi,\chi\ra$.

\subsection{Theorem}
A $G$--Frobenius algebra is naturally a left $k[G]$--module
algebra as well as a right $k[G]$--co--module algebra. Moreover it
satisfies the Yetter--Drinfeld (YD) condition for bi--modules and
is thus a module over $D(k[G])$, the Drinfeld double of $k[G]$.
Where the YD condition reads
\begin{equation}
\sum h_1 \cdot m_0 \otimes h_2m_1 = \sum (h_2 \cdot m)_0\otimes
(h_2\cdot m)_1 h_1
\end{equation}

Here we used the usual notation for co--algebras and right co--modules.
I.e.\ if $\Delta: H \rightarrow H \otimes H$ is the co--multiplication of
$H$ then for $h \in H$ we write $\Delta(h) = \sum h_1 \otimes h_2$
and if $\check\rho: M \rightarrow H$ is a right co--module map, then
for $m\in M$ we write $\check\rho(m)= \sum m_0 \otimes m_1$.

{\bf Proof.} The Theorem follows from the collection of facts below and the
general statement that any $H$ bi--module satisfying the
YD--condition is a a module over $D(H)$ (see e.g.\ [M]).

\subsubsection{Remark} In this particular case the YD condition
states that the co-module structure is $k[G]$ equivariant with
respect to the adjoint action of $k[G]$ on itself, viz.\ as a
tensor product of $G$--Frobenius algebras of left $k[G]$ modules.
See below.

\subsubsection{The $k[G]$--module structure}
Since $A$ is a $k$ algebra,
the $G$--action $\varphi$ turns $A$ into a right $k[G]$ module.
More precisely for $a \in A \sum_g \nu_g g \in k[G]$
\begin{equation}
(\sum \nu_g g) \otimes a \mapsto \sum_g \nu_g \varphi_g(a)
\end{equation}
Since $\varphi \in \Hom(G,Aut(A))$ this is a module structure.

\subsubsection{The $k[G]$--co-module structure}
Since $A$ is a $G$ graded algebra it is naturally a $k[G]$
--co-module.

More precisely for $a \in A, a= \oplus_g a_g$
 the
$k[G]$ co-module structure $\rho:A \rightarrow A \otimes k[G]$ is
given by

\begin{equation}
a \mapsto \sum_g (a_g \otimes g)
\end{equation}

which obviously yields a co-module.

\subsubsection{Lemma} A $G$--Frobenius algebra is a $k[G]$--module algebra
and a $k[G]$ co-module algebra or equivalently a $k[G]^*$--module
algebra.

{\bf Proof.} For the module algebra structure notice that:
\begin{itemize}
\item[1)] $A$ is a left $k[G]$ module as noticed before.
\item[2)] The $k[G]$--action induced by $\varphi$ is by definition
by algebra automorphisms, and $\Delta(g) = g \otimes g$ thus
$$
\varphi_g(ab) = \varphi_g(a) \varphi_g(b)
$$
\item[3)] Since the unit is invariant:
$$
\varphi_g(1) = 1 = \eps(g) 1
$$
\end{itemize}

The structure of co-module algebra follows from the fact that
$$
\varphi(a_g b_h) = a_g b_h \otimes gh
$$
which, as is well known, is nothing but the condition of $A$ being
a $G$ graded algebra

\begin{equation}
\label{kgmod} A_g A_h \subset A_{gh}
\end{equation}

\subsubsection{Remark}
Notice that the condition (\ref{kgmod}) is usually given by a
strict inclusion, so that it is usually not $k[G]$--Galois --
which is equivalent to $A_gA_h=A_{gh}$ (cf.\ [M]). In
case it is, the structure of the algebra is particularly
transparent. We will come back to this later.

\subsubsection{The compatibility}
We will view  $k[G]$ as a left  $k[G]$ module using the adjoint action.
Then $A \otimes k[G]$ turns into a left
$k[G]$--module by using the diagonal action. This is the natural left
$k[G]$-module structure on the tensor product of left Hopf
modules.

$$
 (\sum_h \mu_h h )(\sum_g a_g \otimes \nu g)=
\sum_{h,g} \mu_h \nu \varphi(h)(a_g) \otimes ghg^{-1}
$$

\subsubsection{Lemma}
The co-module structure is $k[G]$--equivariant and thus the
co-module map is a map of left $k[G]$ modules where we use the left
adjoint action of $k[G]$ on itself as the left $k[G]$ action.

\begin{eqnarray}
\rho ((\sum_h \mu_h h ) (a))&=&
\rho (\sum_h \mu_h \varphi_h(a))\nn\\
&=&\sum_{h,g} \mu_h \varphi_h(a_g) \otimes hgh^{-1} = (\sum_h
\mu_h h )\cdot(\sum_g a_g \otimes g)\nn\\
&=&(\sum_h \mu_h h )\cdot \rho(a)
\end{eqnarray}

\subsubsection{The YD condition}
Plugging in the co--product and action yields
\begin{equation}
\varphi_g(a_h) \otimes gh = \varphi_g(a_h) \otimes (ghg^{-1})g
\end{equation}
which verifies the YD condition for $A$.

\subsubsection{Proposition} If $A$ is a  $G$--Frobenius algebra that
is $k[G]$--Galois over $A_e$ as a k[G]--comodule algebra
then $A$ is special and $\g\in Z^2(G,A^*)$
where $A^*$ are the units of $A$. So in particular $\g$ determines $\varphi$
uniquely.

Moreover if $A_e$ is one--dimensional  then $A = k^{\a}[G]$, for some
$\a \in H^2(G,k^*)$ with a choice of parity $\tilde{} \in
\Hom(G,\Z2)$.

{\bf Proof.} Since $A_{g^{-1}}A_g= A_e$ there are elements
$a_g \in A_g, b_{g^{-1}}\in A_{g^{-1}}$ s.t.\ $b_{g^{-1}}a_g=1$
then $a_g$ is a cyclic generator since $\forall c_g \in A_G$
$c_g = c_g(b_{g^{-1}}a_g) = ( c_gb_{g^{-1}}) a_g$ and
$ c_g b_{g^{-1}}\in A_e$. Choosing generators $1_g$ in this way it
is easy to check that the cocycles need to be invertible
and thus the $\varphi$ are fixed by (\ref{grpcompat}). Furthermore
notice that the multiplication map induces an
isomorphism of $A_e$ modules between $A_e$ and $A_g$ via
$a \mapsto a1_g$ where $A_e$ is a cyclic $A_e$ module over itself via left
multiplication. This follows by associativity from $a=a(1_g1_{g^{-1}})
=(a1_g)1_{g^{-1}}$ and thus $a1_g \neq 0$ and the map
$A_e \rightarrow A_g$  is also injective.
Thus the restriction maps are all isomorphisms and
graded cocycles coincide with the usual ones.

\section{The action of Discrete Torsion}
\subsection{Twisting $G$--Frobenius algebras}

Given a $G$--Frobenius algebra
$A$ we can re-scale the multiplication and $G$--action
by a scalar. More precisely let $\l:G\times G \rightarrow k^*$ be
a function. For $a = \oplus _g a_g \in A$ we define

$$\varphi^{\l}(g)(a) = \oplus_h \l(g,h) \varphi(g)(a_h)$$

Given another function $\mu:G\times G \rightarrow k^*$
 we can also define a new multiplication $\circ^{\mu}$
$$
a_g \circ^{\mu} b_h = \mu(g,h) a_g \circ b_h
$$

\subsubsection{Remark}
These twists arise from a projectivization of the
$G$--structures induced on a module over $A$ as for instance the
associated Ramond--space (cf.\ [K1]). In physics terms this means
that each twisted sector will have a projective vacuum, so that fixing their
lifts in different ways induces the twist. Mathematically this means that
$g$ twisted sector is considered to be a Verma module over $A_g$ based on
this vacuum.

\subsubsection{Induced shift on the metric} Due to the
invariance of the metric,
the twist in the multiplication results in a twisted metric
$$\eta^{\mu}(a_g, b_g^{-1}) := \mu(g,g^{-1})\eta (a_g, b_g^{-1})$$

\subsubsection{Definition} We define
$s(\mu,\l)(A)$ to be the collection\\
$\la G,A,\circ^{\mu},1,\eta^{\l},\varphi^{\l},\chi\ra$.

\subsubsection{Proposition}
$s(\mu,\l)(A)$ is
$G$--Frobenius algebra if and only if the following equations
hold for $\mu,\l$:

\begin{equation}
\label{unit}
\mu(e,g)=\mu(g,e)=1
\end{equation}

furthermore $\forall g,h,k \in G$ s.t.\ $A_gA_hA_k \neq 0$:

\begin{equation}
\label{ass}
\mu(g,h)\mu(gh,k)=\mu(h,k)\mu(g,hk)
\end{equation}
and if $A_g A_h \neq 0$ then

\begin{equation}
\label{lmu} \l(g,h) =\frac{\mu(g,h)}{\mu(ghg^{-1},g)}
\end{equation}

If $A_gA_h \neq 0$ as well as $A_gA_hA_k=0$
\begin{equation}
\label{equi}
\l(g,hk)\mu(h,k)=\l(g,h)\l(g,k)\mu(ghg^{-1},gkg^{-1})
\end{equation}

Furthermore
\begin{eqnarray}
\label{transform}
\l(gh,k)& =& \l(h,k)\l(g,hkh^{-1})\nn\\
\l(e,g)&=&\l(g,g)=1\nn\\
\mu([g,h],hgh^{-1})\l(h,g)&=& \mu([g,h],h)\l(g^{-1},ghg^{-1})
\end{eqnarray}
where the third equation has to hold for all pairs $g,h$ s.t.\
$\exists c \in A_{[g,h]}$ s.t.\ $\chi_{h}\mathrm {Tr} (l_c
\varphi_{h}|_{A_{g}})\neq 0$, where $l_c$ is the left multiplication by
$c$. In particular it must hold for all pairs $g,h$ with $[g,h]=e$.

{\bf Proof.}
The first equation (\ref{unit}) expresses that $1$ is still
the unit for the algebra.
The statement (\ref{ass}) for $\mu$ is the obvious form of
associativity. The statement
(\ref{lmu}) comes from the compatibility equation of the group action
with the multiplication.

The equation (\ref{equi}) ensures the equivariance
of the multiplication. It is automatic if $A_gA_hA_k \neq 0$ and
also if $A_g A_h=0$.

The first
equation  (\ref{transform}) for $\l$ is equivalent to the
fact that $\varphi^{\l}$ is
still a $G$ action.

 Notice that $\l(e,g) =1$ since
$A_eA_g=A_g$ and thus
$$
\l(e,g) = \frac{\mu(e,g)}{\mu(g,e)}=1
$$
and so the identity remains
invariant.

Also notice that there is no twist to the character!
$$
\chi^{\l}_g =(-1)^{\tilde g} \dim A_g \Str^{-1}(\varphi_g^{\l)}|_{A_e})
= \chi_g \l(g,e)= \chi_g
$$
This in turn implies the second statement in the second line by projective
self--invariance

$$\chi_{g}^{-1}id|A_{g}=
\varphi^{\l}_{g}|A_{g}=\l(g,g)\varphi_{g}|A_{g}=\l(g,g)\chi_{g}^{-1}id|A_{g}
$$

$$1=\l(e,k)=\l(g^{-1}g,k)=\l(g,k)\l(g^{-1},gkg^{-1})$$
so
$$\l(g,h)=
\l(g^{-1},ghg^{-1})^{-1}$$

The third equation  follows from the projective trace axiom.

$\forall c \in A_{[g,h]}$ and $l_c$ left multiplication by $c$:

\begin{equation}
\label{true} \chi_{h}\mathrm {Tr} (l_c  \varphi_{h}|_{A_{g}})=
\chi_{g^{-1}}\mathrm  {Tr}(  \varphi_{g^{-1}} l_c|_{A_{h}})
\end{equation}

Thus we must have $\chi_{h}\mathrm {Tr} (l_c
\varphi_{h}^{\l}|_{A_{g}})= \chi_{g^{-1}}\mathrm  {Tr}(
\varphi_{g^{-1}}^{\l} l_c|_{A_{h}})$ but this is equivalent to the
third equation in view of (\ref{true}).

Now we check the other axioms.

The invariance of the metric follows from associativity.
\begin{eqnarray*}
\eta^{\mu}(a_g,b_h\circ^{\mu}c_{h^{-1}g^{-1}})&=&
\mu(g,g^{-1})\mu(h, h^{-1}g^{-1})\eta(a_g,b_h\circ c_{h^{-1}g^{-1}})\\
&=&
\mu(g,h)\mu(gh, h^{-1}g^{-1})\eta(a_g \circ b_h,c_{h^{-1}g^{-1}})\\
&=&\eta^{\mu}(a_g \circ^{\mu} b_h,c_{h^{-1}g^{-1}})
\end{eqnarray*}

The projective invariance of the metric  reads as
$$
\l(g,k)\l(g,k^{-1})\mu(gkg^{-1},gk^{-1}g^{-1})=
\mu(k,k^{-1})
$$
which is automatic since $A_g A_k A_{k^{-1}}\neq 0$.


\subsubsection{Definition} We call a twist $s(\l,\mu)$ {\em universal}
if it transforms {\em any}
$G$--Frobenius algebra into a $G$--Frobenius algebra.
We call two twists $s(\l,\mu)$ and $s(\l',\mu')$ isomorphic
if for any $G$--Frobenius algebra $A$
the algebras $s(\l,\mu)(A)$ and $s(\l',\mu')(A)$ are isomorphic.

\subsubsection{Theorem} The universal twists are in 1--1 correspondence with
elements $\a \in Z^2(G,k^*)$  and the isomorphism classes of universal
twists are given by $H^2(G,k^*)$.

{\bf Proof.} If the twist is universal then there are no restrictions on
the equations. In particular $\mu \in Z^2(G,k^*)$ and $\l$ is completely
determined by $\mu$ via (\ref{lmu}). All the other properties are then
automatic. The  claim about isomorphism classes is obvious by noticing
that if $\a$ and $\a'$ are cohomologous and $\a/\a'=d\b$
for some $\b\in Z^1(G,k^*)$ then
a diagonal rescaling of the generators of $k^{\a}[G]$
by $\b$ yields $k^{\a'}[G]$ so the result follows from the charcterization
of universal twists as taking tensor produce with twisted group rings below.

\subsection{Discrete torsion}
In this subsection we prove that universal twists are exactly given by
twisting with discrete torsion.

\subsubsection{Reminder} Given two $G$--Frobenius algebras
$\la G,A,\circ,1,\eta,\varphi,\chi\ra$ and
$\la G,A',\circ',1',\eta',\varphi',\chi'\ra$ we
defined [K1] their tensor product
as $G$--Frobenius algebras to be the $G$--Frobenius algebra

$\la G,\bigoplus_{g \in G}( A_g \otimes A'_g),
\circ\otimes \circ',1\otimes 1',\eta\otimes \eta',\varphi\otimes \varphi',
\chi\otimes \chi'\ra$.

We will use the short hand notation $A \hat \otimes B$ for this product.

\subsubsection{Definition} Given a $G$--Frobenius algebra
$A$ and an element $\a \in H^2(G,k)$ we define the
$\a$--twist of $A$ to be the $G$--Frobenius algebra
$A^{\a}:= A \hat\otimes k^{\a}[G]$.

Notice that
\begin{equation}
\label{alphaiso}
A^{\a}_{g}= A_g \otimes k \simeq A_g
\end{equation}
Using this identification the $G$--Frobenius structures are given by:

\subsection{Lemma} The induced structures under the isomorphism
(\ref{alphaiso}) are
\begin{eqnarray}
\circ^{\a}|_{A^{\a}_{g}\otimes A^{\a}_{h}}&=& \a(g,h) \circ\nn\\
\varphi^{\a}_g|_{A^{\a}_h}&=&\eps(g,h)\varphi_g\nn\\
\eta^{\a}|_{A^{\a}_g\otimes A^{\a}_{g^{-1}}}&=& \a(g,g^{-1})\eta\nn\\
\chi_g=\chi_g
\end{eqnarray}

{\bf Proof.} We notice that the two algebras have the same linear
structure $A_{\a,g}=A_g \otimes kg \simeq A_g$ with
the isomorphism given by $a_g\otimes g \mapsto a_g$. Now the multiplication
is given by
$$(a_g \otimes g) \otimes (a_h \otimes h)
\mapsto a_g a_h\otimes\a(g,h)gh
=  \a(g,h) a_g a_h\otimes gh$$ which yields the twisted
multiplication.

The twist for the $G$ action is computed to be
$$
\varphi_{\a,h}(a_g \otimes g) = \eps(g,h) \varphi_h(g)\otimes
hgh^{-1}
$$

This leads us to the Proposition
\subsubsection{Proposition}
$A_{\a} \simeq s(\a,\eps)A$.

\subsubsection{Theorem}
\label{dt}The set of universal twists
 are described by tensoring with twisted group algebras
which identifies this operation with twisting by discrete torsion.

In other words given a generic $G$--Frobenius algebra $A$ there are
exactly $H^2(G,k)$ twists of it by discrete torsion.

\subsubsection{Discrete torsion as phases for the partition sum}
Notice that  for any
$c\in A^{\a}_{[g,h]}\simeq A_{[g,h]}$

\begin{equation}
\chi_h \Str(l_c \varphi_h|_{A_g^{\a}})=
\eps(h,g) \chi_h \Str(l_c \varphi_h|_{A_g})
\end{equation}
This is the original freedom of choice of a phase for the summands
of the partition function postulated by physicists. In this context,
we should regard $g,h: [g,h]=e$ and set $c=e$. More precisely
set

\begin{equation}
Z(A) = \sum_{g,h \in G :[g,h]=e}  \Str(\chi_g \varphi_g|_{A_h^{\a}})
:= \sum_{g,h \in G :[g,h]=e} Z_{g,h}
\end{equation}
\begin{equation}
Z(A_{\a}) =\sum_{g,h \in G :[g,h]=e} \eps(g,h) Z_{g,h}
\end{equation}

We could omit the factors $\chi_g$, but from the point of view of
physics we should take the trace in the Ramond space (cf.\ [K1]) where
the $k[G]$ module structure is twisted by $\chi$.

\subsection{Supergrading}
\label{super}

In this subsection we wish to address questions of super--grading.
There is a general theory of super--graded $G$--Frobenius algebras
and special $G$--Frobenius algebras. We will expose the structures
here for the group ring.

\subsubsection{Super $G$--Frobenius algebras}
If the underlying algebra of a $G$--Frobenius algebra has a supergrading
$\tilde{}$ then the axioms of a $G$--Frobenius algebra have to be changed to

\begin{itemize}

\item[b$^{\sigma}$)] {\em Twisted super--commutativity}

$a_{g}\circ a_{h} = (-1)^{\tilde a_g\tilde a_h} \varphi_{g}(a_{h})\circ a_{g}$

\item[iv$^{\sigma}$)]
{\em Projective super--trace axiom}

$\forall c \in A_{[g,h]}$ and $l_c$ left multiplication by $c$:

$\chi_{h}\mathrm {STr} (l_c  \varphi_{h}|_{A_{g}})=
\chi_{g^{-1}}\mathrm  {STr}(  \varphi_{g^{-1}} l_c|_{A_{h}})$
\end{itemize}
where $\mathrm{STr}$ is the super--trace.

For details on the super structure as well as
the role of the super structure for special $G$--Frobenius
algebras we refer to [K1].
\subsubsection{Supergraded twisted group rings}
Fix $\a \in H^2(G,k^*), \s \in \Hom(G,\Z2)$ then there is a twisted
super--version of the group ring where now the relations
read
\begin{equation}
 \hat g \hat h =  \a(g,h)\widehat {gh}
\end{equation}
and the twisted commutativity is
\begin{equation}
 \hat g \hat h = (-1)^{\s(g)\s(h)}\varphi_{g}(\hat h) \hat g
\end{equation}
and thus
\begin{equation}
 \varphi_{g}(\hat h)=
(-1)^{\s(g)\s(h)}\a(g,h)\a(gh,g^{-1}) \widehat{ghg^{-1}} =:
\varphi_{g,h} \widehat{ghg^{-1}}
\end{equation}
and thus
\begin{equation}
\eps(g,h) := \varphi_{g,h} = (-1)^{\s(g)\s(h)}\frac{\a(g,h)}{\a(ghg^{-1},g)}
\end{equation}

 We would just like to remark that the
axiom iv$^{\sigma}$ shows the difference between super twists and
discrete torsion.

\subsubsection{Definition} We denote the $\a$-twisted group ring
with super--structure $\s$ by $k^{\a,\s}[G]$.
We still denote $k^{\a,0}[G]$ by $k^{\a}[G]$
where $0$ is the zero map and we denote $k^{0,\s}[G]$
just by $k^{\s}[G]$ where $0$ is the unit of the group $H^2(G,k^*)$.

A straightforward calculation shows
\subsubsection{Lemma} $k^{\a,\s}[G] = k^{\a}[G]\otimes k^{\s}[G]$
and more generally
\subsubsection{Lemma} Let $A$ be
the $G$--Frobenius algebra or more generally super Frobenius algebra with
super grading $\tilde{}$
$\la G,A,\circ,1,\eta,\varphi,\chi\ra$ then
$A\otimes k^{\s}[G]$ is isomorphic to the super $G$--Frobenius algebra
$\la G,A,\circ,1,\eta,\varphi^{\s},\chi^{\s}\ra$ with super grading
$\tilde \s$, where
\begin{equation}
\varphi^{\s}_{g,h} = (-1)^{\s(g)\s(h)}\varphi_{g,h},
\quad \chi^{\s} = (-1)^{\s(g)}\chi_g,
\quad \tilde a_g^{\s}= \tilde a_g + \s(g)
\end{equation}

Using arguments and definitions for universal
twists as for discrete torsion we can obtain the
following Proposition. Here universal means that there is no assumption
on the particular structure of the $G$--Frobenius algebra, in other words
it pertains to generic $G$--Frobenius algebras.
\subsubsection{Proposition} Given a (super) $G$--Frobenius algebra $A$
the universal super $G$--Frobenius algebra re--gradings are
in 1-1 correspondence with $\Hom(G,\Z2)$ and these
structures can be realized by tensoring with $k^{\s}[G]$
for $\s \in  \Hom(G,\Z2)$.

\section{Projective representations, Extensions and twisted group algebras}
 In this section we first assemble classical
facts about groups which will be extended to $G$--Frobenius algebras.
As an intermediate step we analyze twisted group algebras, which belong
to both worlds.

\subsection{Part I: Groups}

\subsubsection{Projective representations}
A projective representation $\rho$ of a group is a map
$\rho:G \rightarrow GL(V)$, $V$ being a $k$--vector space,
which satisfies

\begin{equation}
\rho(g)\rho(h)=\a(g,h) \rho(gh), \rho(e)= id
\end{equation}

It is easy to check that $\a(g,h) \in Z^2(G,k^*)$.
Moreover with
a natural notion of projective isomorphy two projective representations
are isomorphic if their classes are cohomologous (cf.\ e.g.\ [CR, Kar]).

\subsubsection{Extensions}
Given a central extension

\begin{equation}
1 \mapsto A \rightarrow G^* \stackrel{\pi}{\rightarrow} G \rightarrow 1
\end{equation}
fix a section $s$ of $\pi$ and define $\a:G\times G \rightarrow A$ by
$s(g)s(h)= \a(g,h)s(gh)$. It is easy to see that
indeed $\a \in Z^2(G,A)$ and furthermore changing the section or
changing the extension by an isomorphism preserves
the cohomology class of $\a$.

Vice versa a cycle in $\a \in Z^2(G,A)$ were $A$ is an Abelian group
gives rise to a central group extension of $G$.

\begin{equation}
1 \mapsto A \rightarrow G^{\a} \stackrel{\pi}{\rightarrow} G \rightarrow 1
\end{equation}
where $G^{\a}= A\rtimes G$.
The maps are given by $a \mapsto (A,e_G)$ , $(a,g) \mapsto g$ and
the multiplication   is given by
$(a,g) (a',g')= (aa'\a(g,g'),gg')$.

\subsubsection{The transgression map} Given a cycle $\a \in H^2(G,A)$
there is a natural map
\begin{equation}
Tra_{\a}:  \Hom(A,k^*) \rightarrow H^2(G,k^*)
\end{equation}
which sends $\chi\in \Hom(A,k^*)$ to the cocycles defined by
$(g,h)\mapsto \chi\a(g,h)$. Actually this map maps into the cohomology group with values in the torsion subgroup of $k^*$ which we call $tors(k^*)$.
\begin{equation}
Tra_{\a}:  \Hom(A,k^*) \rightarrow H^2(G,tors(k^*))
\end{equation}
\subsubsection{Facts}
\label{fact} We briefly give the facts linking
group cohomology, projective representations and twisted group algebras.
For a detailed account see [Kar].
\begin{itemize}
\item[1)] The classes of central extensions of a group $G$ by an Abelian group
$A$ are in 1--1 correspondence with $H^2(G,A)$.
\item[2)] Any projective $\a$--representation  is a
module over the $\a$--twisted group algebra $k^{\a}[G]$.
(This is in fact an equivalence
of categories).
\item[3)] Every projective representation with
cycle $\a$ is projectively equivalent to one that
can be lifted to linear representation on $G^{\hat\a}$
if $[\a]$ is in the image of the transgression map associated to $[\hat \a]$.
\item[4)] If $H^2(G,k^*)= H^2(G,(tors(k^*))$ then:
\begin{itemize}
\item[a)]  any projective representation can be lifted to a suitable group, and
\item[b)] there is a universal extension
$$
1 \rightarrow A \rightarrow G^* \rightarrow G \rightarrow 1
$$
such that {\em any} projective representation lifts to $G^*$ and
moreover the group $A\simeq H^2(G,k^*)$.
\end{itemize}
\end{itemize}

\medskip

\noindent {\sc Assumption:}\\
\nopagebreak
For the remainder of the section we will assume that  $k$
has the property that $H^2(G,k^*)= H^2(G,(tors(k^*))$. This
is the case e.g.\ if is algebraically
closed or $k= \mathbb{R}$, see eg.\ [Kar].

\subsection{Part II: The twisted group algebra revisited}
\label{partii}

Fix $[\a']\in H^2(G,A)$, an element $[\a]\in \mathrm{Im}(Tra_{[\a']})$
and a pre--image character $\chi \in \Hom(A,k^*)$.

This yields a central extension:
\begin{equation}
1 \mapsto A \rightarrow G^{\a'} \stackrel{\pi}{\rightarrow} G \rightarrow 1
\end{equation}
with a section $s$ of $\pi$ s.t.\ the cocycle corresponding to
$s$ is $\a'$. The map  $\chi$ induces map
\begin{equation}
\chi: k[G^{\a}]\rightarrow k[G] : ag\mapsto \chi(a)g
\end{equation}

while the section $s$ induces a map

\begin{equation}
s: k[G]\rightarrow k[G^{\a}]:  g \mapsto 1_Ag
\end{equation}

\subsubsection{Projective algebra}
Using the maps $s,\chi$ we can also characterize the multiplication
$\mu^{\a}$ in $k^{\a}[G]$
as follows: It is the map which makes the following diagram commutative:
$$
\begin{CD}
k[G^{\a}] \otimes  k[G^{\a}] @>\mu>> k[G^{\a}]\\
@A s \otimes s AA @VV \chi V\\
k[G] \otimes  k[G] @>\mu^{\a}>> k[G]\\
\end{CD}
$$

We already know that $\mu^{\a}$ induces the structure of an algebra.
This diagram captures the statement about lifts of projective representations
of $G$ to linear representations of $G^{\a}$.

This is essentially \ref{fact} 1.

\subsubsection{Projective co--algebra}
Using the diagram above as
we define a co--multiplication by commutativity of:

$$
\begin{CD}
k[G^{\a}] @>\Delta>> k[G^{\a}] \otimes  k[G^{\a}]\\
@A sAA @VV \chi\otimes \chi V\\
k[G] @>\Delta^{\a}>> k[G]\otimes  k[G] \\
\end{CD}
$$
The co--algebra structure we induce in this way on $k[G]$ is actually
the old co--algebra structure, but $k^{\a}[G]$ ceases to be a
bi--algebra.

\subsubsection{Remark: Braiding} If one would like a bi--algebra
structure on the group ring $k^{\a}[G]$ then one has to consider
braided objects, where the braiding is inverse to the twist.
It should be possible to find analogous statements to the ones presented
in this article by considering structures over $k^{\a}[G]$
in braided categories.

\subsubsection{Adjoint action}
Let $ad$ denote the adjoint action $k[G^{\a}]$. Then there is an induced
action on $k[G]$
$$
\begin{CD}
k[G^{\a}] \otimes  k[G^{\a}] @>ad>> k[G^{\a}]\\
@A s \otimes s AA @VV \chi V\\
k[G] \otimes  k[G] @>ad^{\eps}>> k[G]\\
\end{CD}
$$

According to \ref{adjointaction} this action is given by

$$ad^{\eps}(g)(h):= \eps(g,h) ghg^{-1}$$

\subsection{Part III: $G$--Frobenius algebras}

We now apply the logic of part II to general $G$--Frobenius algebras.

Let $H$ be an Abelian group.
Fix $[\a']\in H^2(G,H)$, an element $[\a]\in \mathrm{Im}(Tra_{[\a']})$ and a pre-image
character $\chi \in \Hom(H,k^*)$ and a  central extension:
\begin{equation}
1 \mapsto H \rightarrow G^{\a'} \stackrel{\pi}{\rightarrow} G \rightarrow 1
\end{equation}
with a section $s$ of $\pi$ s.t.\ the cocycle corresponding to
$s$ is $\a'$.

\subsubsection{Definition}
Let  $A^{\a}$ be a $G^{\a}$
Frobenius algebra.
We say that a $G$--Frobenius algebra $F$ can be lifted to  $A^{\a}$
if there are maps $i: A \rightarrow A^{\a}$ and
$res: A^{\a} \rightarrow A$
such that the structural maps fit into the commutative
diagrams

$$
\begin{CD}
A^{\a} @>\rho^{\a}>> A^{\a} \otimes  k[G^{\a}]\\
@A i AA @VV res \otimes \chi V\\
A @>\rho>> A \otimes  k[G] \\
\end{CD}
$$

and

$$
\begin{CD}
A^{\a} \otimes A^{\a} @>\mu^{\a}>> A^{\a}\\
@A i\otimes i AA @VV res V\\
A \otimes A@>\mu>> A  \\
\end{CD}
\qquad
\begin{CD}
k[G^{\a}] \otimes A^{\a} @>\varphi^{\a}>> A^{\a}\\
@A s\otimes s AA @VV res V\\
k[G] \otimes A@>\varphi>> A  \\
\end{CD}
$$

and all algebraic structures are compatible.


\subsubsection{Definition} We say that an $H$--Frobenius
algebra $B$ is $H$ homogeneous if it is
endowed  with an additional left $H$--action $\tau$ which shifts group degree
and is equivariant w.r.t.\ multiplication. More precisely the following two
equations hold:
\begin{equation}
\tau(h)(A_{h'}) \subset A_{h,h'}, \qquad \tau(h)(ab)= a \tau(h)(b)
\end{equation}

It is standard to see that

\subsubsection{Remark} With the notation as above,
the left action $\tau$ of $H$ on $B$ is necessarily
by isomorphisms and thus $B$ is a special
$H$--Frobenius algebra whose components are all isomorphic.
Moreover $B$ is Galois as a $k[H]$--comodule over $B_e$.

\subsubsection{Definition} Given a $G$--Frobenius algebra $A$,
an $H$--homogeneous $H$--Frobenius
algebra $B$
and a cocycle $\a\in Z^2(G,H)$ we define the {\em crossed product}
of $A$ and $B$ to be the $G^{\a}$--Frobenius algebra
\begin{equation}
A\#_{\a}B:=\la G^{\a},A\otimes B, \circ \#_{\a}\circ',
1\otimes 1, \eta \otimes \eta', \varphi  \#_{\eps} \varphi', \chi\otimes \chi' \ra
\end{equation}

Where
\begin{equation}
(a_g \otimes b_h) \circ \#_{\a}\circ'(c_{g'} \otimes d_{h'}) =
a_g c_{g'} \otimes \tau(\a(g,g')) b_h d_{h'}
\end{equation}
and
\begin{equation}
 \varphi  \#_{\eps} \varphi'(g,h) (a_{g'} \otimes b_{h'})=
\varphi_g(a_{g'})\otimes \t(\a(g,g')\a(gg',g^{-1}))\varphi_{h}(b_{h'})
\end{equation}

We leave it to the reader to verify all axioms, since it is analogous
to previous calculations.

\subsubsection{Quantum symmetry group} The postulated second left
action by translation $\t$
can be viewed as the quantum symmetry group postulated
by physicists. Notice that it acts freely. The invariants are
linearly isomorphic to $\bigoplus A_{g} \otimes B_{e_H}$ where
$e_H$ is the unit element of $H$.

\subsubsection{Lemma} The linear map above
induces an isomorphism
$$
{}^{H_{\tau}}(A \#_{\a} B) \simeq\bigoplus_{g\in G} (A_{g}  \otimes B_{e_H})
$$
as $G$--Frobenius algebras with trivial action on the second factor.

Here we denoted the invariants under the action of $H$ by
$\t$ by ${}^{H_{\tau}}$.

\subsubsection{Definition}  Fix $\chi\in \Hom(H,k^*)$
then there is a natural map from $B$ to $B_e$ given by
$\t(h)b_e \mapsto \chi(h)b_e$. This map induces a map

$$A \#_{\a} B \rightarrow
 \bigoplus_{g\in G} (A_{g}  \otimes B_{e_H})
$$
which induces a structure of $G$--Frobenius algebra on
$\bigoplus_{g\in G} (A_{g}  \otimes B_{e_H})$, the $G$--action
on the second factor being trivial.
We define
$$(A \#_{\a} B)^{\chi}$$ to be this $G$--Frobenius algebra.
\medskip

It is easy to check that the following holds:
\subsubsection{Lemma} Keeping the notation above, let
$[\a']= Tra_{[\a]}(\chi)$ and more precisely on the level of cocycles
let $\a'(g,g') = \chi\a(g,g')$.

Then

$$
(A \#_{\a} B)^{\chi} \simeq \left(\bigoplus (A_{g}  \otimes B_{e_H})\right)_{\a'}
$$

\subsubsection{Definition} Given a cocycle  $\a\in Z^2(G,H)$, a
central extension $G^{\a}$ of $G$ by $H$
and a $G$--Frobenius algebra $A$ we define $A^{\a}$ to be the
$G^{\a}$--Frobenius algebra
$$
A^{\a}:= A\#_{\a} k[H]
$$

\subsubsection{Theorem}

Given a $G$--Frobenius algebra $A$ and  cocycles $\a \in  Z^2(G,k^*), \a' \in Z^2(G,H)$
which are related by $\chi \in \Hom(G,k^*)$ via
$\a(g,g') = \chi(\a'(g,g'))$.

Then the twist $A_{\a}$ of $A$ lifts to the $G^{\a'}$--Frobenius algebra
$A^{\a'}$
and moreover
$$
(A^{\a})^{\chi} \simeq A_{\a}
$$

Finally if $G^*$ is the universal extension of $G$ whose
cocycle is $\b \in H^2(G, H^2(G,k^*))$ then any twist $A_{\a}$
of a  $G$--Frobenius
algebra $A$  lifts to $A^{\b}$.

{\bf Proof.} Choose a section $s$ of the
extension yielding $\a$. We denote the unit element
of $H$ by $e_H$ and denote $s(g)$ by $e_Hg$.
We let $i: A_g \rightarrow A^{\a'}_{e_Hg}$ be the map given
by $A_g \rightarrow A_g\otimes k e_H: a_g \mapsto a_g \otimes e_H$
and define $res: A^{\a'}_{hg}\simeq
A_g \otimes k e_H \mapsto A_g$ to be the map
 $a_g \otimes e_H \mapsto \chi(h) a_g$.

Then

\begin{multline*}
(res\otimes \chi)( \rho^{\a'} (i(a_g)))=
(res\otimes \chi)( \rho^{\a'} (a_g \otimes e_H))\\
=(res\otimes \chi)(a_g \otimes e_H) \otimes (e_Hg)
= a_g \otimes g
\end{multline*}
which assures the comodule algebra structure.

$$
\chi ( \mu^{\a} ((i\otimes i)(a_g \otimes b_{g'})))
= \chi(a_{g} b_{g'} \otimes \a(g,g')gg')=
\a(g,h)a_g b_{g'}
$$
since $A_{e_Hg}A_{e_Hg'}\subset A_{\a(g,g')gg'}$
which assures the algebra structure.

$$
\chi \circ \varphi \circ (s\otimes s)(g\otimes a_h)
=\chi \circ \varphi ((e_Hg\otimes A_{e_Hg'}))
= \eps(g,g')\varphi_{g}(a_{g'})
$$
which assures the module algebra structure,
since $\varphi_{e_Hg}(A_{e_Hh}) \subset A_{\eps'(g,h)gh}$,
where we set
$$
\eps'(g,g') =\frac{\a'(g,h)}{\a'(ghg^{-1},g)}
$$
to be the cocycle of the adjoint action.
Then by  \ref{adjointaction}
$$
\chi(\eps'(g,g'))=\eps(g,g')
$$

For the last statement notice that (cf.\ e.g.\ [Kar])

$$k[G^*] = \prod_{\a\in T} k[G^{\a}]$$
where $T$ is a transversal for $B^2(G,k^*)$ in $Z^2(G,k^*)$.
So $A^{\b} \simeq \bigoplus_{\a\in T} A^{\a}$
and we can lift to the appropriate component.


\begin{thebibliography}{99}
\bibitem
 [A]  {A}
P.\ A.\ Aspinwall
 {\it A Note on the Equivalence of Vafa's and
Douglas's Picture of Discrete Torsion}
JHEP 0012 (2000) 029.

\bibitem
 [CR]  {CR}
C.\ W.\ Curtis and
I.\ Reiner.
{\it Methods of representation theory. Vol. I.
 With applications to finite groups and orders.}
 Pure and Applied Mathematics.  John Wiley \& Sons, Inc., New York, 1988.
\bibitem
 [D]  {D}
M.\ R.\ Douglas.
{\it  D-branes and Discrete Torsion}
Preprint hep-th/9807235.

\bibitem
 [FG]  {FG}

 B.\ Fantechi and L.\ Goettsche.
 {\it Orbifold cohomology for global quotients.} Preprint. math.AG/0104207


\bibitem
 [Kar]  {Kar}
G.\ Karpilovsky.
{\it Projective Representations of Finite Groups}. Dekker, 1985.

\bibitem
 [K1]  {K1}
R.\ Kaufmann. {\it Orbifolding Frobenius algebras}.
Preprint. MPI 2001-56, IHES M/01/37, math.AG/0107163 and
 R.\ Kaufmann. {\it Orbifolding Frobenius algebras}.
Talk at WAGP2000 conference at SISSA Trieste, October 2000

\bibitem
 [K2]  {K2}
R.\ Kaufmann. {\it Second quantized Frobenius algebras}.
Preprint. IHES M/02/46, math.AG/020613.
\bibitem
 [K3]  {K3}
R.\ Kaufmann. {\it Singularities with symmetries, orbifold Frobenius algebras
and mirror symmetry.}
(in preparation)
\bibitem
 [LS]  {LS}

M.\ Lehn and C.\ Soerger.
{\it The cup product of the Hilbert scheme for K3 surfaces}
Preprint math.AG/0012166.



\bibitem
 [M]  {M}
S.\ Montgomery.
{\it  Hopf algebras and their actions on rings}.
CBMS Regional Conference Series in Mathematics, 82.
American Mathematical Society, Providence, RI, 1993.
\bibitem
 [R]  {R}
Y.\ Ruan
{\it Discrete torsion and twisted orbifold cohomology}.
Preprint math.AG/0005299

\bibitem
 [S]  {S}
E.\ Sharpe.
{\it Discrete Torsion and Gerbes I+II}
  Preprints hep-th/9909120 and hep-th/9909108

\bibitem
 [V]  {V}
C.\ Vafa
{\it Modular Invariance and Discrete Torsion on Orbifolds}.
 Nucl. Phys. B273 (1986) 592--606.
\bibitem
 [VW]  {VW}
C.\ Vafa and E.\ Witten.
{\it On orbifolds with discrete torsion.}
 J.\ Geom.\ Phys.\ 15 (1995), 189--214.
\end{thebibliography}
\end{document}